\newif\iffigs\figstrue
\documentclass[12pt]{article}
\iffigs
   \input{epsf}
\else
\fi

\newcommand{\sect}[1]{\setcounter{equation}{0}\section{#1}}

\newcommand{\eq}{\begin{equation}}
\newcommand{\eqa}{\begin{eqnarray}}
\newcommand{\en}{\end{equation}}
\newcommand{\ena}{\end{eqnarray}}
\newcommand{\enn}{\nonumber \end{equation}}

\def\sk{\vskip .4cm}
\def\noi{\noindent}
\def\om{\omega}

\def\al{\alpha}
\def\la{\lambda}

\def\ga{\gamma}
\def\Ga{\Gamma}

\def\Cb{\bar{C}}

\def\rhop{{\rho}^{\prime}}

\def\epsi{\varepsilon}
\def\we{\wedge}

\def\de{\delta}

\def\part{\partial}

\def\R#1#2{ R^{#1}_{~~~#2} }

\def\L#1#2{ \La^{#1}_{~~~#2} }

\def\La{\Lambda}

\def\Cb{{\bf \mbox{\boldmath $C$}}}

\def\c#1#2{ C^{#1}_{~~#2} }
\def\C#1#2{ {\bf \mbox{\boldmath $C$}}^{#1}_{~~#2} }

\def\limepsizero{\lim_{\epsi \rightarrow 0}}

\def\n2{{{N+1} \over 2}}

\def\square{{\,\lower0.9pt\vbox{\hrule \hbox{\vrule height 0.2 cm
\hskip 0.2 cm \vrule height 0.2 cm}\hrule}\,}}

\def\Q.E.D.{\rightline{$\Box$}}

\def\sumong{\sum_{g \in G}}

\def\sumongnote{\sum_{g \not= e}}

\def\sumonh{\sum_{h \in G}}
\def\sumonhp{\sum_{h \in G'}}

\def\Lcal{{\cal L}}
\def\Rcal{{\cal R}}

\def\omc#1#2{\om^{#1}_{~~#2}}
\def\W#1#2{W^{#1}_{~~~#2} }

\def\Gc#1#2{\Ga^{#1}_{~~#2}}
\def\a#1#2{a^{#1}_{~#2}}
\def\ainv#1#2{(a^{-1})^{#1}_{~#2}}

\begin{document}
\begin{titlepage}
\vskip -1cm
\rightline{DFTT-19/2000}
%\rightline{LBNL-40824}
%\rightline{August 1999}
\vskip 1em
\begin{center}
{\large\bf Differential calculi on finite groups }
\\[2em]
Leonardo Castellani\\[.4em]
{\sl Dipartimento di Scienze e Tecnologie Avanzate,
 East Piedmont University, Italy; \\
Dipartimento di Fisica Teorica and Istituto Nazionale di
Fisica Nucleare\\
Via P. Giuria 1, 10125 Torino, Italy.} \\[2em]
\end{center}
\vskip 4 cm
\begin{abstract}
 A brief review of bicovariant differential calculi on finite groups
is given, with some new developments on diffeomorphisms and
integration. We illustrate the general theory
with the example of the nonabelian finite group $S_3$.
\end{abstract}
%\noi{April 1998}~~\hskip 10cm~q-alg/9709032

\vskip 8cm
\centerline{{\small publ. in JHEP, Proceedings of the}}
\centerline{{\small {\bf Corfu
Summer Institute on Elementary Particle Physics, 1998}}}
\vskip 1cm
\noi \hrule
\vskip .2cm
\noi{\small e-mail: castellani@to.infn.it}
%\noi$^*${\small \sl II Facolt\`a di Scienze M.F.N., sede di Alessandria}
\end{titlepage}
\newpage
\setcounter{page}{1}
%%%%%%%%%%%%%%%%%%%%%%%%%%%%%%
\sect{Introduction}

Differential calculi can be constructed on spaces
that are more general than differentiable manifolds.
Indeed
the algebraic construction of differential calculus in terms of Hopf
structures allows to extend the usual differential geometric
quantities (connection, curvature, metric, vielbein etc.) to a
variety of interesting spaces that include quantum groups,
noncommutative spacetimes (i.e. quantum cosets), and discrete
spaces.
 \sk
 In this contribution  we concentrate on the differential
 geometry of finite group ``manifolds". As we will
discuss, these spaces can be visualized as collections of points,
corresponding to the finite group elements, and connected by
oriented links according to the particular differential calculus
we build on them. Although functions $f \in Fun(G)$ on finite
groups $G$ commute, the calculi that are constructed on $Fun(G)$
by algebraic means are in general noncommutative, in the sense
that differentials do not commute with functions, and the exterior
product does not coincide with the usual antisymmetrization of the
tensor product.
 \sk
 Among the physical motivations for finding
differential calculi on finite groups we mention
the possibility of using finite group spaces
as internal spaces for Kaluza-Klein compactifications of Yang-Mills,
(super)gravity or superstring theories (
for example Connes'
reconstruction of the standard model in terms of noncommutative
geometry \cite{Connes}  can be recovered as Kaluza-Klein
compactification of Yang-Mills theory on an appropriate discrete
internal space). Differential calculi on discrete
spaces can be of use in the study of integrable models, see for ex. ref.
\cite{DMintegrable}. Finally gauge
and gravity theories on finite group spaces may be used as lattice
approximations. For example the action for pure Yang-Mills $\int F
\we {}^* F$ considered on the finite group space $Z^N \times Z^N
\times Z^N \times Z^N$, yields the usual Wilson action of lattice
gauge theories, and $N \rightarrow \infty$ gives the continuum
limit \cite{DMgauge}. New lattice theories can be found by choosing
different finite groups.

 A brief review of the differential calculus
on finite groups is presented. Most of this material is not new,
and draws on the treatment of ref.s \cite{DMGcalculus,FI1,gravfg},
 where the Hopf
algebraic approach of Woronowicz \cite{Wor} for the construction of
differential calculi is adapted to the setting of finite groups.
Some developments on Lie derivative, diffeomorphisms and
integration are new. The general theory is illustrated in the case
of $S_3$.
%%%%%%%%%%%%%%%%%%%%%%%%%%%%%%
\sect{Differential calculus on finite groups}

Let $G$ be a finite group of order $n$ with
generic element $g$ and unit $e$. Consider $Fun(G)$, the set of
complex functions on $G$. An element $f$ of $Fun(G)$ is specified
by its values $f_g \equiv f(g)$ on the group elements $g$, and can
be written as \eq f=\sum_{g \in G} f_g x^g,~~~f_g \in \Cb
\en
where the functions $x^g$ are defined by
\eq
x^g(g') = \de^g_{g'}
\en
Thus $Fun(G)$ is a n-dimensional vector space, and the $n$ functions
$x^g$ provide a basis. $Fun(G)$ is also a commutative algebra,
with the usual pointwise sum and product [$(f+h)(g)=f(g)+h(g)$,
$(f\cdot h)(g)=f(g)h(g)$, $ (\la f)(g)=\la f(g), f,h \in Fun(G),
\la \in \Cb$] and unit $I$ defined
by $I(g)=1, \forall g \in G$. In particular:
\eq
x^g x^{g'}=\de_{g,g'} x^g,~~~\sumong x^g = I \label{mul}
\en
Consider now the
left multiplication by $g_1$:
 \eq L_{g_1}g_2=g_1g_2,~~~\forall
g_1,g_2 \in G
\en
This induces the left action (pullback) $\Lcal_{g_1}$ on
$Fun(G)$:
\eq
\Lcal_{g_1} f(g_2) \equiv f(g_1g_2)|_{g_2},~~~\Lcal_{g_1}:Fun(G)
\rightarrow Fun(G)
\en
where $f(g_1g_2)|_{g_2}$ means $f(g_1g_2)$ seen as a function
of $g_2$.
Similarly we can define the right action  on $Fun(G)$ as:
 \eq
(\Rcal_{g_1}f)(g_2)= f(g_2g_1)|_{g_2}
\en
For the basis functions we find easily:
\eq
\Lcal_{g_1} x^{g} = x^{g_1^{-1} g},
~~\Rcal_{g_1} x^{g} = x^{g g_1^{-1}}
\en
Moreover:
 \eqa & &\Lcal_{g_1} \Lcal_{g_2}=\Lcal_{g_1g_2},
~~\Rcal_{g_1} \Rcal_{g_2}=\Rcal_{g_2g_1},\\ & &\Lcal_{g_1}
\Rcal_{g_2}=\Rcal_{g_2} \Lcal_{g_1}
 \ena
  \sk
\noi
{\bf Bicovariant differential calculus}
\sk
Differential calculi can be constructed on Hopf algebras $A$ by
algebraic means, using the costructures of $A$ \cite{Wor}.
In the case of finite groups $G$, differential calculi on $A=Fun(G)$ have been
discussed in ref.s \cite{DMGcalculus,FI1,gravfg}. Here we give the main results
derived in  \cite{gravfg}, to which we refer for a more detailed treatment.
\sk

A {\bf {\sl first-order differential calculus}} on $A$ is defined by
\sk
i) a linear map $d$: $A \rightarrow \Gamma$, satisfying the Leibniz rule
\eq
d(ab)=(da)b+a(db),~~\forall a,b\in A; \label{Leibniz}
\en
The ``space of 1-forms" $\Ga$ is an
appropriate bimodule on $A$, which
essentially means that its elements can be
multiplied on the left and on the right by elements of $A$
[more precisely $A$ is a left module if $\forall a,b \in A, \forall
\rho,\rho' \in \Ga $ we have: $ a(\rho+\rho')=a\rho+a\rho',
~(a+b)\rho=a\rho+b\rho, ~a(b\rho)=(ab)\rho,~ I\rho=\rho$. Similarly
one defines a right module. A left and right module
is a {\sl bimodule} if
$a(\rho b)=(a\rho)b$]. From the Leibniz rule
$da=d(Ia)=(dI)a+Ida$ we deduce $dI=0$.
\sk
ii) the possibility of expressing any $\rho \in \Ga$ as
\eq
\rho=\sum_k a_k db_k \label{adb}
\en
\noi for some $a_k,b_k$ belonging to $A$.
\sk
 To build a first
order differential calculus on $Fun(G)$ we need to extend the
algebra $A=Fun(G)$ to a differential algebra of elements
$x^g,dx^g$ (it is sufficient to consider the basis elements and
their differentials). Note however that the $dx^g$ are not
linearly independent. In fact from $0=dI=d(\sumong x^g)=\sumong
dx^g$ we see that only $n-1$ differentials are independent. Every
element $\rho = adb$ of $\Ga$ can be expressed as a linear
combination (with complex coefficients) of terms of the type $x^g
dx^{g'}$. Moreover $\rho b \in \Ga$ (i.e. $\Ga$ is also a right
module) since the Leibniz rule and the multiplication rule
(\ref{mul}) yield the commutations:
 \eq
  dx^g x^{g'} = -x^g
dx^{g'}+\de^g_{g'} dx^g
\en
allowing to reorder functions to the left of differentials.
 \sk
 \noi {\bf Partial derivatives}
 \sk
  Consider the differential of a
function $f \in Fun(g)$: \eq df = \sumong f_g dx^g = \sumongnote
f_g dx^g + f_e dx^e= \sumongnote (f_g - f_e)dx^g \equiv
\sumongnote
\part_g f dx^g \label{partcurved}
\en
We have used $dx^e = - \sumongnote dx^g$ (from $\sumong dx^g=0$).
The partial derivatives of $f$ have been defined in analogy with
the usual differential calculus, and are given by \eq
\part_g f = f_g - f_e = f(g) - f(e) \label{partcurved2}
\en
Not unexpectedly, they take here the form of finite differences
(discrete partial derivatives at the origin $e$).
\sk
\noi
 {\bf Left and right covariance}
\sk
A differential calculus is
 {\sl left or right covariant} if the left or right action of
 $G$ ($\Lcal_g$ or $\Rcal_g$) commutes with the exterior derivative $d$.
 Requiring left and right covariance in fact {\sl defines} the action of
 $\Lcal_g$ and $\Rcal_g$ on differentials: $\Lcal_g db \equiv
 d(\Lcal_g b), \forall b \in Fun(G)$ and similarly for
 $\Rcal_g db$. More generally, on elements of $\Ga$
 (one-forms) we define $\Lcal_g$ as:
 \eq
 \Lcal_g (adb) \equiv (\Lcal_g a) \Lcal_g db =
 (\Lcal_g a) d (\Lcal_g b)
 \en
 and similar for $\Rcal_g$.
 Computing for example the left and right action on the differentials
 $dx^g$ yields:
 \eq
\Lcal_g (dx^{g_1})\equiv d(\Lcal_g x^{g_1})=dx^{g^{-1}g_1},~~
\Rcal_g (dx^{g_1})\equiv d(\Rcal_g x^{g_1})=dx^{g_1 g^{-1}}
\en
A differential calculus is called {\sl bicovariant} if it is
both left and right covariant.
\sk
\noi
{\bf Left invariant one forms}
 \sk
 As in usual Lie group manifolds, we can introduce a basis in $\Ga$
 of left-invariant one-forms $\theta^g$:
 \eq
  \theta^g \equiv
\sumonh x^{hg} dx^h ~~(=\sumonh x^h dx^{hg^{-1}}),
\label{deftheta}
\en
It is immediate to check that $\Lcal_k \theta^g = \theta^g$. The
relations (\ref{deftheta}) can be inverted:
 \eq
  dx^h = \sumong (x^{hg} - x^h)\theta^g \label{dxastheta}
 \en
{} From $0=dI=d\sumong x^g =\sumong dx^g=0$ one finds:
\eq
 \sumong
\theta^g = \sumong \sumonh x^h dx^{hg^{-1}}= \sumonh x^h \sumong
dx^{hg^{-1}}=0 \label{sumtheta}
\en
Therefore we can take as basis of the cotangent space $\Ga$ the
$n-1$ linearly independent left-invariant one-forms $\theta^g$
with $g \not= e$ (but smaller sets of $\theta^g$ can be
consistently chosen as basis, see later).
 \sk
  The commutations
between the basic 1-forms $\theta^g$ and functions $f \in Fun(G)$
are given by: \eq f \theta^g = \theta^g \Rcal_g f
\label{fthetacomm}
\en
Thus functions do commute between themselves (i.e. $Fun(G)$ is
a commutative algebra) but do not commute with the
basis of one-forms $\theta^g$. In this sense the differential geometry
of $Fun(G)$ is noncommutative, the noncommutativity being
milder than in the case of quantum groups $Fun_q(G)$(which
are noncommutative algebras).
\sk
The right
action of $G$ on the elements $\theta^g$ is given by:
\eq
\Rcal_h \theta^g = \theta^{ad(h)g},~~\forall h \in G
\en
where $ad$ is the adjoint action of $G$ on $G$, i.e. $ad(h)g
\equiv hgh^{-1}$. Then {\sl bicovariant calculi are in 1-1
correspondence with unions of conjugacy classes (different from
$\{e\}$)} \cite{DMGcalculus}: if $\theta^g$ is set to zero, one must set to
zero all the $\theta^{ad(h)g},~\forall h \in G$ corresponding to the
whole conjugation class of $g$.
 \sk
 We denote by $G'$ the subset corresponding
 to  the union of conjugacy classes
 that characterizes the bicovariant calculus on $G$
 ($G' = \{g \in G |\theta^g \not= 0\}$).
 Unless otherwise indicated, repeated indices are
 summed on $G'$ in the following.
 \sk

A bi-invariant (i.e. left and right invariant)
one-form $\Theta$ is obtained by summing
on all $\theta^g$ with $g \not= e$:
\eq
\Theta = \sumongnote \theta^g
\en
 {\bf Exterior product}
\sk
For a bicovariant differential calculus on a Hopf
algebra $A$ an {\sl exterior product}, compatible with the left
and right actions of $G$, can be defined by
\eq
\theta^{g_1} \we \theta^{g_2}=\theta^{g_1} \otimes \theta^{g_2}
- \theta^{g_1^{-1} g_2 g_1} \otimes \theta^{g_1} \label{exprod}
\en
where the tensor product between elements $\rho,\rhop \in \Ga$
is defined to
have the properties $\rho a\otimes \rhop=\rho \otimes a \rhop$, $a(\rho
\otimes \rhop)=(a\rho) \otimes \rhop$ and $(\rho \otimes \rhop)a=\rho
\otimes (\rhop a)$.

Note that:
 \eq
  \theta^{g} \we \theta^{g}=0~~~~\mbox{(no sum on $g$)}
\en

  Left and right actions on $\Ga \otimes \Ga$ are
  simply defined by:
  \eq
  \Lcal_h (\rho \otimes \rhop)= \Lcal_h \rho \otimes \Lcal_h
  \rhop,~~~
\Rcal_h (\rho \otimes \rhop)= \Rcal_h \rho \otimes \Rcal_h
  \rhop
  \en
  (with the obvious generalization to $\Ga \otimes ...\otimes \Ga$)
  so that for example:
\eq
\Lcal_h (\theta^i \otimes \theta^j)=
\theta^i \otimes \theta^j,~~~~ \Rcal_h (\theta^i \otimes \theta^j)=
\theta^{ad(h)i} \otimes \theta^{ad(h)j}
\en
  We can generalize the definition
(\ref{extheta})to exterior products of $n$ one-forms:
 \eq
\theta^{i_1} \we ... \we \theta^{i_n} \equiv \W{i_1i_2}{j_1k_1}
\W{k_1 i_3}{j_2k_2} \W{k_2 i_4}{j_3k_3} ...
\W{k_{n-2}i_n}{j_{n-1}j_n} \theta^{j_1} \otimes ...\otimes
\theta^{j_n}
\en
where the matrix $W$ is defined by:
\eq
 \theta^i \we \theta^j \equiv \W{ij}{kl} \theta^k \otimes \theta^l=
 \theta^i \otimes \theta^j -
\L{ij}{kl} \theta^k \otimes \theta^l. \label{extheta}
\en
and $\L{ij}{kl}$ is the braiding matrix defined by
(\ref{exprod}).
The space of $n$-forms $\Ga^{\we n}$ is therefore defined as in
the usual case but with the new permutation operator $\La$,
and can be shown to be a bicovariant bimodule, with left
and right action defined as for $\Ga \otimes ...\otimes \Ga$
with the tensor product replaced by the wedge product.
\sk
\noi
{\bf Exterior derivative}
\sk
 Having the exterior product we can define the {\sl exterior
derivative}
\eq
d~:~\Gamma \rightarrow \Gamma \we \Gamma
\en
\eq
d (a_k db_k) = da_k \we db_k,
\en
\noi which can easily be extended to $\Gamma^{\we n}$ ($d:
\Gamma^{\we n} \rightarrow \Gamma^{\we (n+1)}$), and has the
following properties:
\eq
 d(\rho \we \rhop)=d\rho \we \rhop +
(-1)^k \rho \we d\rhop \label{propd1}
\en
\eq
d(d\rho)=0\label{propd2}
\en
\eq
\Lcal_g (d\rho)=d \Lcal_g \rho \label{propd3}
\en
\eq \Rcal_g (d\rho)=d \Rcal_g \rho \label{propd4}
\en
\noi where $\rho \in \Ga^{\we k}$, $\rhop \in \Ga^{\we n}$. The
last two properties express the fact that $d$ commutes with the
left and right action of $G$.
\sk
\noi
{\bf Tangent vectors}
\sk
Using (\ref{dxastheta}) to expand $df$ on the basis of
the left-invariant one-forms
$\theta^g$ defines the (left-invariant) tangent vectors $t_g$:
 \eq
  df=\sumong f_g dx^g  = \sumonhp (\Rcal_{h^{-1}} f - f ) \theta^h
\equiv \sumonhp (t_h f) \theta^h \label{partflat}
\en
so that the ``flat" partial derivatives $t_h f$ are given by
\eq
t_h f = \Rcal_{h^{-1}} f - f \label{partflat2}
\en
 The Leibniz rule for the flat partial derivatives $t_g$
reads:
\eq
 t_g (ff')=(t_g f) \Rcal_{g^{-1}} f'  + f t_g f' \label{tgLeibniz}
\en

In analogy with ordinary differential calculus, the operators
$t_g$ appearing in (\ref{partflat}) are called (left-invariant)
{\sl tangent vectors}, and in our case are given by
 \eq
  t_g =
\Rcal_{g^{-1}}- id \label{tangent}
\en
They satisfy the composition rule: \eq t_g t_{g'}= \sum_h
\c{h}{g,g'} t_h \label{chichi}
\en
where the structure constants are:
\eq
\c{h}{g,g'}=\de^h_{g'g} - \de^h_{g}-\de^h_{g'}
\label{cconst}
\en
and have the property:
\eq
\c{ad(h)g_1}{~~ad(h)g_2,ad(h)g_3}= \c{g_1}{g_2,g_3} \label{adhc}
\en

{\bf Note 2.1 :}
The exterior derivative on any $f \in Fun(G)$ can be expressed as
a commutator of $f$ with the bi-invariant one-form $\Theta$:
\eq
df = [\Theta , f]
\en
as one proves by using (\ref{fthetacomm}) and (\ref{partflat}).
\sk
 {\bf Note 2.2 :} From the fusion rules (\ref{chichi}) we
deduce the ``deformed Lie algebra" (cf. ref.s \cite{Wor,ACintro,
Athesis}):
 \eq
  t_{g_1} t_{g_2} -
\L{g_3,g_4}{g_1,g_2}t_{g_3} t_{g_4}= \C{h}{g_1,g_2} t_h
\en
where the $\Cb$ structure constants are given by:
 \eq
\C{g}{g_1,g_2} \equiv \c{g}{g_1,g_2} - \L{g_3,g_4}{g_1,g_2}
\c{g}{g_3,g_4}= \c{g}{g_1,g_2} -
\c{g}{g_2,g_2 g_1 g_2^{-1}}=
\de^{ad(g_2^{-1})g}_{g_1} - \de^g_{g_1}
\label{Cconst}
\en
and besides property (\ref{adhc}) they also satisfy:
 \eq
\C{g}{g_1,g_2}=\C{g_1}{g,g_2^{-1}} \label{propC}
\en
Moreover the following identities hold:

{\bf i)} {\sl deformed Jacobi identities:}
\eq \C{k}{h_1,g_1}
\C{h_2}{k,g_2} - \L{g_3,g_4}{g_1,g_2}\C{k}{h_1,g_3}\C{h_2}{k,g_4}=
\C{k}{g_1,g_2} \C{h_2}{h_1,k} \label{Jacobi}
\en

{\bf ii)} {\sl fusion identities:} \eq
\C{k}{h_1,g} \C{h_2}{k,g'}=
\c{h}{g,g'} \C{h_2}{h_1,h} \label{adfusion}
\en

Thus the $\Cb$ structure constants are a representation (the
adjoint representation) of the tangent vectors $t$. \sk \noi {\bf
Cartan-Maurer equations, connection and curvature} \sk From the
definition (\ref{deftheta}) and eq. (\ref{fthetacomm}) we deduce
the Cartan-Maurer equations:
 \eq
 d\theta^g + \sum_{g_1,g_2}
\c{g}{g_1,g_2}\theta^{g_1}\we \theta^{g_2}=0 \label{CM}
\en
where the structure constants $\c{g}{g_1,g_2}$ are those
given in (\ref{cconst}).
\sk
Parallel transport of the vielbein $\theta^g$
can be defined as in
ordinary Lie group manifolds:
\eq
\nabla \theta^g= - \omc{g}{g'} \otimes \theta^{g'}
 \label{parallel}
\en
where $\omc{g_1}{g_2}$ is the connection one-form:
\eq
\omc{g_1}{g_2}= \Gc{g_1}{g_3,g_2} \theta^{g_3}
\en
Thus parallel transport is a map from $\Ga$ to $\Ga \otimes \Ga$;
by definition it must satisfy:
\eq
\nabla (a \rho) = (da)\otimes \rho + a \nabla \rho,~~~\forall a \in
A,~\rho \in \Ga \label{parallel1}
\en
and it is a simple matter to verify that this relation is
satisfied with the usual parallel transport of Riemannian
manifolds. As for the exterior differential, $\nabla$ can be
extended to a map $\nabla : \Ga^{\we n} \otimes \Ga
\longrightarrow \Ga^{\we (n+1)} \otimes \Ga $ by defining:
 \eq
\nabla (\varphi \otimes \rho)=d\varphi \otimes \rho +
 (-1)^n \varphi \nabla
\rho
\en

Requiring parallel transport to commute with the left and right
action of $G$ means:
 \eqa & &\Lcal_{h} (\nabla \theta^{g})=\nabla
( \Lcal_{h} \theta^{g}) =\nabla \theta^g\\ & &\Rcal_{h} (\nabla
\theta^{g})=\nabla ( \Rcal_{h} \theta^{g}) =\nabla \theta^{ad(h)g}
\ena
 Recalling that  $\Lcal_{h} (a \rho)=(\Lcal_h a) (\Lcal_h
\rho)$ and $\Lcal_{h} (\rho \otimes \rho')=(\Lcal_h \rho) \otimes
(\Lcal_h \rho'),~\forall a \in A,~\rho,~\rho' \in \Ga$ (and
similar for $\Rcal_h$),
 and substituting
(\ref{parallel}) yields respectively:
\eq
\Gc{g_1}{g_3,g_2} \in \Cb
\en
and
 \eq
 \Gc{ad(h)g_1}{ad(h)g_3,ad(h)g_2}=\Gc{g_1}{g_3,g_2} \label{adga}
\en
Therefore the same situation arises as in the case of Lie groups,
for which  parallel transport on the group manifold commutes
with left and right action iff the connection components are
$ad(G)$ - conserved constant tensors. As for Lie groups, condition
(\ref{adga}) is satisfied if one takes $\Ga$ proportional to the
structure constants. In our case, we can take any combination of
the $C$ or $\Cb$ structure constants, since both are $ad(G)$
conserved constant tensors. As we see below, the $C$ constants
can be used to define a torsionless connection, while the $\Cb$
constants define a parallelizing connection.

\sk
 As usual, the {\sl curvature} arises from $\nabla^2$:
  \eq
  \nabla^2 \theta^g = - \R{g}{g'} \otimes \theta^{g'}
\en
\eq
\R{g_1}{g_2} \equiv d \omc{g_1}{g_2} + \omc{g_1}{g_3} \we
\omc{g_3}{g_2} \label{curvature}
\en

The {\sl torsion} $R^g$ is defined by:
\eq
R^{g_1} \equiv d\theta^{g_1} +  \omc{g_1}{g_2} \we \theta^{g_2}
\label{torsion}
\en

Using the expression of $\om$ in terms of $\Ga$ and the
Cartan-Maurer equations yields
 \eqa
 \R{g_1}{g_2} &=& (-
\Gc{g_1}{h,g_2} \c{h}{g_3,g_4} + \Gc{g_1}{g_3,h} \Gc{h}{g_4,g_2})~
\theta^{g_3} \we \theta^{g_4}=\\
&=& (-
\Gc{g_1}{h,g_2} \C{h}{g_3,g_4} + \Gc{g_1}{g_3,h} \Gc{h}{g_4,g_2}-
\Gc{g_1}{g_4,h} \Gc{h}{g_4g_3g_4^{-1},g_2})~\theta^{g_3} \otimes
 \theta^{g_4}\nonumber
\ena
\eqa
R^{g_1}&=& (- \c{g_1}{g_2,g_3} + \Gc{g_1}{g_2,g_3})~ \theta^{g_2} \we
\theta^{g_3}= \nonumber \\
&=& (- \C{g_1}{g_2,g_3} + \Gc{g_1}{g_2,g_3}-
\Gc{g_1}{g_3,g_3g_2g_3^{-1}})~\theta^{g_2} \otimes
\theta^{g_3}
\ena

Thus a connection satisfying:
 \eq
  \Gc{g_1}{g_2,g_3}-
\Gc{g_1}{g_3,g_3g_2g_3^{-1}}=\C{g_1}{g_2,g_3} \label{rconn}
  \en
corresponds to a vanishing torsion $R^g =0$ and could be
  referred to as a ``Riemannian" connection.
\sk
 On the other hand,  the choice:
   \eq
  \Gc{g_1}{g_2,g_3}=\C{g_1}{g_3,g_2^{-1}} \label{parconn}
  \en
corresponds to a vanishing curvature $\R{g}{g'}=0$, as can be
checked by using the fusion equations (\ref{adfusion}) and
property (\ref{propC}). Then (\ref{parconn}) can be called the
parallelizing connection: {\sl finite groups are
parallelizable.}
\sk
\noi
{\bf Tensor transformations }
\sk
Under the familiar transformation of the connection 1-form:
\eq
(\omc{i}{j})' = \a{i}{k} \omc{k}{l} \ainv{l}{j} +
\a{i}{k} d \ainv{k}{j} \label{omtransf}
\en
the curvature 2-form transforms homogeneously:
\eq
(\R{i}{j})' = \a{i}{k} \R{k}{l} \ainv{l}{j}
\en
The transformation rule (\ref{omtransf}) can be seen as induced by
the change of basis $\theta^i=\a{i}{j} \theta^j$, with $\a{i}{j}$
invertible $x$-dependent matrix (use eq. (\ref{parallel1}) with
$a\rho=\a{i}{j} \theta^j$).
 \sk
  \noi
   {\bf Metric}
\sk
 The metric tensor $\ga$ can be defined as an element of $\Ga
\otimes \Ga$:
 \eq
  \ga = \ga_{i,j} \theta^i \otimes \theta^j
  \en
 Requiring it to be invariant under left and right action of
 $G$ means:
 \eq
 \Lcal_h (\ga)=\ga=\Rcal_h (\ga)
 \en
or equivalently, by recalling $\Lcal_h(\theta^i \otimes
\theta^j)=\theta^i \otimes \theta^j$, $\Rcal_h(\theta^i \otimes
\theta^j)=\theta^{ad(h)i}\otimes \theta^{ad(h)j}$  :
 \eq
  \ga_{i,j} \in \Cb,~~  \ga_{ad(h)i,ad(h)j}=\ga_{i,j} \label{gabiinv}
 \en
These properties are  analogous to the ones satisfied by the
Killing metric of Lie groups, which is indeed constant and
invariant under the adjoint action of the Lie group.
 \sk
 On finite $G$ there are various choices of biinvariant
 metrics. One can simply take $\ga_{i,j}=\de_{i,j}$,
   or $\ga_{i,j}= \C{k}{l,i} \C{l}{k,j}$.
   \sk
 For any biinvariant metric $\ga_{ij}$ there are tensor transformations
 $\a{i}{j}$ under which $\ga_{ij}$ is invariant, i.e.:
 \eq
 \a{h}{h'} \ga_{h,k} \a{k}{k'}=\ga_{h',k'} \Leftrightarrow
\a{h}{h'} \ga_{h,k} = \ga_{h',k'} \ainv{k'}{k} \label{ginv}
 \en
 These transformations are simply given by the matrices that
rotate the indices according to the adjoint action of $G$:
 \eq
 \a{h}{h'} (g) = \de^{ad(\al(g))h}_{h'} \label{Gadjoint}
 \en
 where $\al(g): G \mapsto G$ is an arbitrary mapping. Then
 these matrices are functions of $G$ via this mapping, and
 their action leaves $\ga$ invariant because of the its biinvariance
 (\ref{gabiinv}). Indeed
substituting these matrices in (\ref{ginv}) yields:
 \eq
 \a{h}{h'} (g) \ga_{h,k} \a{k}{k'} (g)=
 \ga_{ad([\al(g)]^{-1})h',ad([\al(g)]^{-1})k'}= \ga_{h',k'}
 \en
proving the invariance of $\ga$.

Consider now a contravariant vector $\varphi^i$ transforming as
$(\varphi^i)'=\a{i}{j}(\varphi^j)$. Then using (\ref{ginv}) one
can easily see that
 \eq
  (\varphi^k \ga_{k,i})'=  \varphi^{k'} \ga_{k',i'} \ainv{i'}{i}
  \en
  i.e.  the vector $\varphi_i \equiv \varphi^k \ga_{k,i}$ indeed
  transforms as a covariant vector.
\sk
 \noi
  {\bf Lie derivative and diffeomorphisms}
   \sk
    The notion of diffeomorphisms, or general coordinate transformations,
is fundamental in gravity theories. Is there such a notion in the
setting of differential calculi on Hopf algebras ? The answer is
affirmative, and has been discussed in detail in ref.s
\cite{ACintro,LCqISO,Athesis}.
As for differentiable manifolds, it relies on the existence of the
Lie derivative.

Let us review the situation for Lie group manifolds. The Lie
derivative $l_{t_i}$ along a left-invariant tangent vector $t_i$
is related to the infinitesimal right translations generated by
$t_i$:
 \eq
 l_{t_i} \rho = \limepsizero {1\over \epsi} [\Rcal_{\exp [\epsi
 t_i]} \rho - \rho] \label{Liederivative1}
 \en
 $\rho$ being an arbitrary tensor field. Introducing the
 coordinate dependence
 \eq
 l_{t_i} \rho (y) = \limepsizero {1\over \epsi}
  [\rho (y + \epsi t_i) - \rho (y) ]
 \en
identifies the Lie derivative $ l_{t_i}$ as a directional
derivative along $t_i$.  Note the difference in meaning of the
symbol $t_i$ in the r.h.s. of these two equations: a group
generator in the first, and the corresponding tangent vector in
the second.
\sk
 For finite groups the Lie derivative takes the form:
 \eq
 l_{t_g} \rho =  [\Rcal_{g^{-1}} \rho - \rho] \label{Liederivative2}
 \en
so that the Lie derivative is simply given by
 \eq
l_{t_g}=\Rcal_{g^{-1}}-id=t_g
\en
cf. the definition of $t_g$ in (\ref{tangent}). For example
 \eq
l_{t_g} (\theta^{g_1} \otimes \theta^{g_2}) =
 \theta^{ad(g^{-1})g_1} \otimes \theta^{ad(g^{-1})g_2}-
 \theta^{g_1} \otimes \theta^{g_2}
 \en

 As in the case of differentiable manifolds, the Cartan formula
 for the Lie derivative acting on p-forms holds:
\eq
 l_{t_g}= i_{t_g} d + d i_{t_g}
 \en
see ref.s  \cite{ACintro,Athesis,gravfg}.

 Exploiting this formula, diffeomorphisms
  (Lie derivatives) along generic tangent vectors $V$
 can also be consistently defined via the operator:
\eq
 l_{V}= i_{V} d + d i_{V}
 \en
 This requires
  a suitable definition
 of the contraction operator $i_V$  along generic tangent vectors
 $V$, discussed in ref. \cite{Athesis,gravfg}.

 We have then a way
 of defining ``diffeomorphisms" along arbitrary (and x-dependent)
 tangent vectors for any tensor $\rho$:
 \eq
 \delta \rho = l_V \rho
 \en
and of testing the invariance of candidate lagrangians under the
generalized Lie derivative.
\sk
\noi
{\bf Haar measure and integration}
\sk
Since we want to be able to define actions (integrals on
$p$-forms) we must now define integration of $p$-forms on finite
groups.

 Let us start with integration of functions $f$. We define the integral
  map $h$ as a linear functional $h: Fun(G) \mapsto \Cb$ satisfying the
  left and right invariance conditions:
  \eq
  h(\Lcal_g f)=0=h(\Rcal_g f)
  \en
  Then this map is uniquely determined (up to a normalization constant),
  and is simply given by the ``sum over $G$" rule:
  \eq
  h(f)= \sumong f(g)
  \en

 Next we turn to define the integral of a p-form.
Within the differential calculus we have a basis of left-invariant
1-forms, which may allow the definition of a biinvariant volume
element. In general for a differential calculus with $n$
independent tangent vectors, there is an integer $p  \geq n$ such
that the linear space of $p$-forms is 1-dimensional, and $(p+1)$-
forms vanish identically. We will see explicit examples in the
next Section. This means that every product of $p$ basis
one-forms $\theta^{g_1} \we \theta^{g_2} \we ... \we \theta^{g_p}$
is proportional to one of these products, that can be chosen to
define the volume form $vol$:
 \eq
 \theta^{g_1} \we \theta^{g_2} \we ... \we \theta^{g_p}=
 \epsilon^{g_1,g_2,...g_p} vol
 \en
 where $\epsilon^{g_1,g_2,...g_p}$ is the proportionality constant.
 Note that the volume $p$-form is obviously left invariant. We can
  prove that it is also right invariant with the following
  argument. Suppose that $vol$ be given by
  $\theta^{h_1} \we \theta^{h_2} \we ... \we \theta^{h_p}$ where
  $h_1,h_2,...h_p$ are given group element labels. Then the right
  action on $vol$ yields:
  \eq
  \Rcal_g [\theta^{h_1} \we  ... \we
  \theta^{h_p}]=
  \theta^{ad(g)h_1} \we ... \we
  \theta^{ad(g)h_p}=
  \epsilon^{ad(g)h_1,...ad(g)h_p} vol
  \en
  Recall now that the ``epsilon tensor" $\epsilon$ is necessarily
  made out of the  $W$ tensors of eq. (\ref{extheta}), defining the wedge
  product. These tensors are invariant under the adjoint action
  $ad(g)$, and so is the $\epsilon$ tensor. Therefore
  $\epsilon^{ad(g)h_1,...ad(g)h_p}=\epsilon^{h_1,...h_p}=1$
  and $\Rcal_g vol = vol$. This will be verified in the examples
  of next Section.

Having identified the volume $p$-form it is natural to set
 \eq
 \int f vol \equiv h(f) = \sumong f(g) \label{intpform}
 \en
 and  define the integral on a $p$-form $\rho$ as:
 \eq
  \int \rho = \int \rho_{g_1,...g_p}~ \theta^{g_1}
 \we ... \we \theta^{g_p}=   \int
\rho_{g_1,...g_p}~\epsilon^{g_1,...g_p} vol \equiv \sumong
\rho_{g_1,...g_p}(g)~\epsilon^{g_1,...g_p}
  \en
Due to the biinvariance of the volume form, the integral map $\int
: \Ga^{\we p} \mapsto \Cb$ satisfies the biinvariance conditions:
 \eq
  \int \Lcal_g f = \int f = \int \Rcal_g f
  \en

  Moreover, under the assumption that the volume form belongs to
  a nontrivial cohomology class, that is $d ( vol) = 0$ but
  $vol \not= d \rho$, the important property holds:
  \eq
  \int df =0
  \en
  with $f$  any $(p-1)$-form: $f=f_{g_2,...g_p}~ \theta^{g_2}
\we ... \we \theta^{g_p}$. This property, which allows
  integration by parts, has a simple proof. Rewrite
  $\int df$ as:
  \eq
\int df= \int (d f_{g_2,...g_p})\theta^{g_2} \we ... \we
\theta^{g_p}+ \int  f_{g_2,...g_p} d ( \theta^{g_2} \we ... \we
\theta^{g_p})
  \en
 Under the cohomology assumption the second term in the r.h.s.
 vanishes, since $d ( \theta^{g_2} \we ... \we
\theta^{g_p}) =0$ (otherwise, being a $p$-form, it should be
proportional to $vol$, and this would contradict the assumption
 $vol \not= d \rho$). Using now (\ref{partflat}) and
 (\ref{intpform}):
\eqa
 & &\int df= \int (t_{g_1} f_{g_2,...g_p})\theta^{g_1}\we \theta^{g_2} \we ...
\we \theta^{g_p} = \int
[\Rcal_{g_1^{-1}}f_{g_2,...g_p}-f_{g_2,...g_p}]
\epsilon^{g_1,...g_p} vol = \nonumber \\
 & &~~~~~~~ = \epsilon^{g_1,...g_p}
 \sumong [\Rcal_{g_1^{-1}}f_{g_2,...g_p}(g)-f_{g_2,...g_p}(g)]=0
\ena
 Q.E.D.

\sect{Bicovariant calculus on $S_3$}

In this Section we illustrate the general theory on
the particular example of the permutation group $S_3$.
 \sk
 Elements: $a=(12)$, $b=(23)$, $c=(13)$, $ab=(132)$, $ba=(123)$, $e$.
 \sk
 Nontrivial conjugation classes: $I = [a,b,c]$, $II = [ab,ba]$.
 \sk
There are 3 bicovariant calculi $BC_I$, $BC_{II}$, $BC_{I+II}$
corresponding to the possible unions of the conjugation classes
\cite{DMGcalculus}. They have respectively dimension 3, 2 and 5.
We examine here the $BC_I$ and $BC_{II}$ calculi. \sk {\bf $BC_I$
differential calculus}
 \sk
 \noi Basis of the 3-dimensional vector space
of one-forms:
 \eq
\theta^a,~\theta^b,~\theta^c
\en
\noi Basis of the 4-dimensional vector space of two-forms:
 \eq
  \theta^a \we \theta^b,~
 \theta^b \we \theta^c,~\theta^a \we \theta^c,~\theta^c \we \theta^b
 \en

 Every wedge product of two $\theta$ can be expressed
 as linear combination of the basis elements:
 \eq
  \theta^b \we \theta^a = -\theta^a \we \theta^c - \theta^c \we
  \theta^b,~~\theta^c \we \theta^a=-\theta^a \we \theta^b-\theta^b \we
  \theta^c
  \en

\noi Basis of the 3-dimensional vector space of three-forms:
 \eq
  \theta^a \we \theta^b \we \theta^c,~\theta^a \we \theta^c \we \theta^b,
  ~\theta^b \we \theta^a \we \theta^c
  \en

 and we have:
  \eqa
 & &\theta^c \we \theta^b \we \theta^a=- \theta^c \we \theta^a \we
\theta^c=-\theta^a \we \theta^c \we \theta^a=\theta^a \we \theta^b
\we \theta^c \nonumber \\
 & &\theta^b \we \theta^c \we \theta^a=- \theta^b \we \theta^a \we
\theta^b=-\theta^a \we \theta^b \we \theta^a=\theta^a \we \theta^c
\we \theta^b \nonumber \\
 & &\theta^c \we \theta^a \we \theta^b=- \theta^c \we
\theta^b \we \theta^c=-\theta^b \we \theta^c \we \theta^b=\theta^b
\we \theta^a \we \theta^c
\ena

\noi Basis of the 1-dimensional vector space of four-forms:

\eq
 vol = \theta^a \we \theta^b \we \theta^a \we \theta^c
 \en

 and we have:
  \eq
 \theta^{g_1} \we \theta^{g_2} \we \theta^{g_3} \we \theta^{g_4}=
 \epsilon^{g_1,g_2,g_3,g_4} vol \label{epsiI}
 \en
 where the nonvanishing components of the $\epsilon$ tensor are:
 \eqa
 & &
 \epsilon_{abac}=\epsilon_{acab}=\epsilon_{cbca}=\epsilon_{cacb}=
 \epsilon_{babc}=\epsilon_{bcba}=1 \\
 & &
 \epsilon_{baca}=\epsilon_{caba}=\epsilon_{abcb}=\epsilon_{cbab}=
 \epsilon_{acbc}=\epsilon_{bcac}=-1 \label{epsivaluesI}
 \ena

\noi Cartan-Maurer equations:

\eqa
 & & d\theta^a+\theta^b\we\theta^c+\theta^c \we \theta^b =0
 \nonumber \\
 & & d\theta^b+\theta^a\we\theta^c+\theta^c \we \theta^a =0
 \nonumber \\
 & & d\theta^c+\theta^a\we\theta^b+\theta^b \we \theta^a =0
 \ena

 The exterior derivative on any three-form
of the type $\theta \we \theta \we \theta$ vanishes, as one can
easily check by using the Cartan-Maurer equations and the
equalities between exterior products given above. Then, as shown
in the previous Section,  integration of a total differential
vanishes on the ``group manifold" of $S_3$ corresponding to the
$BC_I$ bicovariant calculus. This ``group manifold" has three
independent directions, associated to the cotangent basis
$\theta^a,~\theta^b,~\theta^c$. Note however that the volume
element is of order four in the left-invariant one-forms $\theta$.
\sk
{\bf $BC_{II}$ differential calculus}
\sk
\noi Basis of the 2-dimensional vector space of one-forms:
 \eq
\theta^{ab},~\theta^{ba}
\en
\noi Basis of the 1-dimensional vector space of two-forms:
 \eq
  vol = \theta^{ab} \we \theta^{ba}=-\theta^{ba} \we \theta^{ab}
 \en
so that:
  \eq
 \theta^{g_1} \we \theta^{g_2}=
 \epsilon^{g_1,g_2} vol \label{epsiII}
 \en
 where the $\epsilon$ tensor is the usual 2-dimensional Levi-Civita tensor.
 \sk
\noi Cartan-Maurer equations:

\eq
 d\theta^{ab} =0,~~ d\theta^{ba} =0
 \en

 Thus the exterior derivative on any one-form $\theta^g$
 vanishes and integration of a total differential
vanishes on the group manifold of $S_3$ corresponding to the
$BC_{II}$ bicovariant calculus. This group manifold has two
independent directions, associated to the cotangent basis
$\theta^{ab},~\theta^{ba}$.
\sk
{\bf Visualization of the $S_3$ group ``manifold"}
\sk

 We can draw a picture of the group manifold of $S_3$. It is made
 out of 6 points, corresponding to the group elements and identified with
 the functions $x^e,x^a,x^b,x^c,x^{ab},x^{ba}$.
\sk
$BC_I$ - calculus:
\sk
 From each
 of the six points $x^g$ one can move in three directions, associated to
 the tangent vectors $t_a,t_b,t_c$, reaching three other points
 whose ``coordinates" are
 \eq
\Rcal_a x^g = x^{ga},~~\Rcal_b x^g = x^{gb},~~\Rcal_c x^g = x^{gc}
\en
 The 6 points and the ``moves" along the 3 directions are
 illustrated in the Fig. 1. The links are not oriented since
 the three group elements $a,b,c$ coincide with their inverses.

 \sk
$BC_{II}$ - calculus:
\sk
 From each
 of the six points $x^g$ one can move in two directions, associated to
 the tangent vectors $t_{ab},t_{ba}$, reaching two other points
 whose ``coordinates" are
 \eq
\Rcal_{ab} x^g = x^{gba},~~\Rcal_{ba} x^g = x^{gab}
\en
 The 6 points and the ``moves" along the 3 directions are
 illustrated in Fig. 1. The arrow convention on a link
 labeled  (in italic) by a group element $h$
 is as follows: one
 moves in the direction of the arrow via the action
 of $\Rcal_{h}$ on $x^g$. (In this case $h=ab$). To move in the opposite
 direction just take the inverse of $h$.
\sk

 The pictures in Fig. 1
 characterize the bicovariant calculi  $BC_I$ and $BC_{II}$ on $S_3$,
 and were drawn in  ref. \cite{DMGcalculus} as examples of
 digraphs, used to characterize different calculi on sets.
 Here we emphasize  their geometrical meaning as
  finite group ``manifolds".

%\sk \iffigs
%\begin{figure}[h]
%\epsfxsize = 10cm \centerline{\epsfig{figure=S3.eps}}
%\caption{$S_3$ group manifold, and moves of the points under the
%group action}
%\label{s3}
% \unitlength=1mm
%\end{figure}
%\fi

%%Begin InstantTeX Picture
\let\picnaturalsize=N
\def\picsize{4.0in}
\def\picfilename{S3.eps}
%If you do not have the picture file add:
%\let\nopictures=Y
%to the beginning of the file.
\ifx\nopictures Y\else{\ifx\epsfloaded Y\else\input epsf \fi
\let\epsfloaded=Y
\centerline{\ifx\picnaturalsize N\epsfxsize \picsize\fi
\epsfbox{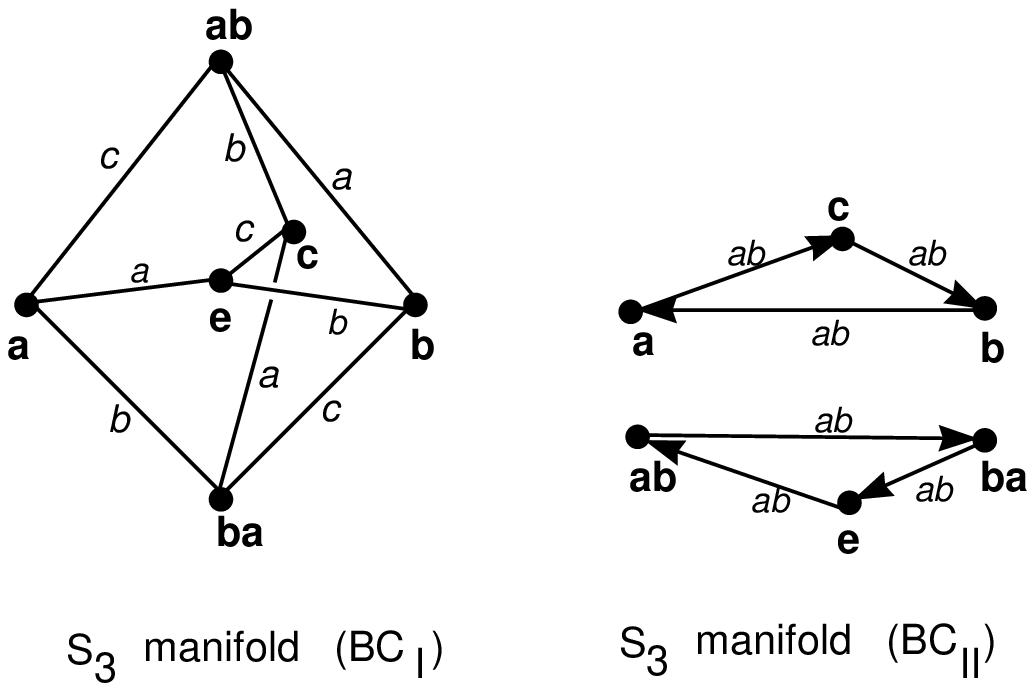}}}\fi
%%End InstantTeX Picture
\sk
  {\small{\bf Fig. 1} : $S_3$ group manifold, and moves of the
   points under the group action}
 \sk
{\bf Acknowledgements }
\sk
It is a pleasure to thank the organizers of the
 Corfu Summer Institute on  Elementary Particle Physics
 for their invitation to discuss physics in such
 a beautiful and relaxed atmosphere.

%%%%%%%%%%%%%%%%%%%%%%%%%%%%%%%%%%%

\vfill\eject
\end{document}